\documentclass[11pt,english]{article}
\usepackage{times}
\usepackage[T1]{fontenc}
\usepackage{a4wide}
\usepackage{graphicx}
\usepackage{setspace}
\usepackage{amssymb}

\makeatletter

\providecommand{\LyX}{L\kern-.1667em\lower.25em\hbox{Y}\kern-.125emX\@}

 \usepackage{verbatim}

\date{1st June 2006}

\usepackage{cite}

\newtheorem{theorem}{Theorem}
\newtheorem{lemma}{Lemma}

\newcounter{mpnkt}
\newcommand{\punkt}{\refstepcounter{mpnkt} \arabic{mpnkt}.~}

\newcommand{\ds}{\displaystyle}

\newcommand{\esli}{\mbox{ \ if \ }}

\usepackage{babel}
\makeatother
\begin{document}

\newcommand{\nN}{\mathcal{N}_{N}}

\newcommand{\ioik}{\left(i_{1},\ldots ,i_{k}\right)}

\newcommand{\Ls}{L_{\mbox {\scriptsize s}}\mbox {}}

\newcommand{\Lsii}{L_{\mbox {\scriptsize s},\ioik }}

\newcommand{\RR}{\mathbb{R}}

\newcommand{\sE}{\mathsf{E}}

\newcommand{\al}{\alpha }

\newcommand{\de}{\delta }

\newcommand{\mux}{\boldmath {\mu }}

\newcommand{\erx}{R}

\newcommand{\Lf}{L_{0}\, }

\newcommand{\cdv}{\varkappa }

\newcommand{\om}{\omega }

\newcommand{\gn}{\gamma _{N}}

\newcommand{\mud}{\mu }

\newcommand{\dsp}{d}

\newcommand{\sdv}{b_{2}}

\newcommand{\srM}{s}

\title{Long-time behavior of stochastic model with multi-particle synchronization }

\author{Anatoly~Manita%
\footnote{This work is supported by Russian Foundation of Basic Research (grants
06-01-00662 and 05-01-22001).  \protect \\
Address:   Faculty of Mathematics and Mechanics, Moscow State University,
119992, Moscow, Russia. \protect \\
E-mail:~manita@mech.math.msu.su%
}}

\maketitle
\begin{abstract}
We consider a basic stochastic particle system consisting of $N$
identical particles with isotropic $k$-particle synchronization,
$k\geq 2$. In the limit when both number of particles $N$ and time
$t=t(N)$ grow to infinity we study an asymptotic behavior of a coordinate
spread of the particle system. We describe three time stages of $t(N)$
for which a qualitative behavior of the system is completely different.
Moreover, we discuss the case when a spread of the initial configuration
depends on~$N$ and increases to infinity as~$N\rightarrow \infty $.
\end{abstract}

\paragraph*{\punkt Introduction.}

A wide class of probabilistic models can be interpreted as stochastic
particle systems with a synchronization-like interaction. In general
terms, the matter concerns a special class of jump Markov processes
$x(t)=\left(x_{1}(t),\ldots ,x_{N}(t)\right)$ evolving in continuous
time and taking their values in $\RR ^{N}$, with generators of the
following symbolic form $L=\Lf +\Ls \, .$ The variable $x_{i}\in \RR $
is interpreted as a coordinate of the particle~$i$. A free dynamics
$\Lf $ corresponds to independent movements of the individual particles,
and $\Ls $ corresponds to synchronizing jumps $x=(x_{1},\ldots ,x_{N})\rightarrow x'=(x'_{1},\ldots ,x'_{N})$
which, according to many preceding papers~\cite{Man-Sche,malysh-Manita-tvp,malysh-Manita-v,Malyshkin},
fit to the following general rule\[
\left\{ x'_{1},\ldots ,x'_{N}\right\} \subset \left\{ x_{1},\ldots ,x_{N}\right\} ,\quad \left\{ x'_{1},\ldots ,x'_{N}\right\} \not =\left\{ x_{1},\ldots ,x_{N}\right\} .\]

First mathematical papers on stochastic synchronization dealt with
the case $N=2$ (see, for example, \cite{Mitra-Mitr}). Many-particle
systems $(N>2)$, studied till now, have different forms of the synchronizing
interaction. So~\cite{Man-Sche,malysh-Manita-tvp,malysh-Manita-v}
considered \emph{pairwise} interactions, while in~\cite{Malyshkin}
a \emph{three-particle} interaction was studied. In the present paper
we consider a \emph{general $k$-particle} symmetrized interaction
which will be defined in terms of synchronizing maps.

When we analyze a collective behavior of such particle system we face
with superposition of two opposite tendencies: with the course of
time the free dynamics increases the spread of the particle system
while the synchronizing interaction tries to decrease it. The aim
of this paper is to study a qualitative balance between these two
tendencies in the limit $N\rightarrow \infty $, $t=t(N)\rightarrow \infty $.
Our main result (Theorem~\ref{t-R}) determines three phases of different
behavior of the system depending on the growth rate of~$t(N)$. Similar
result was obtained in~\cite{malysh-Manita-tvp} for a system with
\emph{two types} of particles and a \emph{pairwise} interaction between
different types.

Main results of this paper were presented by the author at~PTAP-2006~(see~\cite{Man-petroz}).

\paragraph*{\punkt Synchronizing maps.}

\begin{figure}
\includegraphics[  scale=0.5]{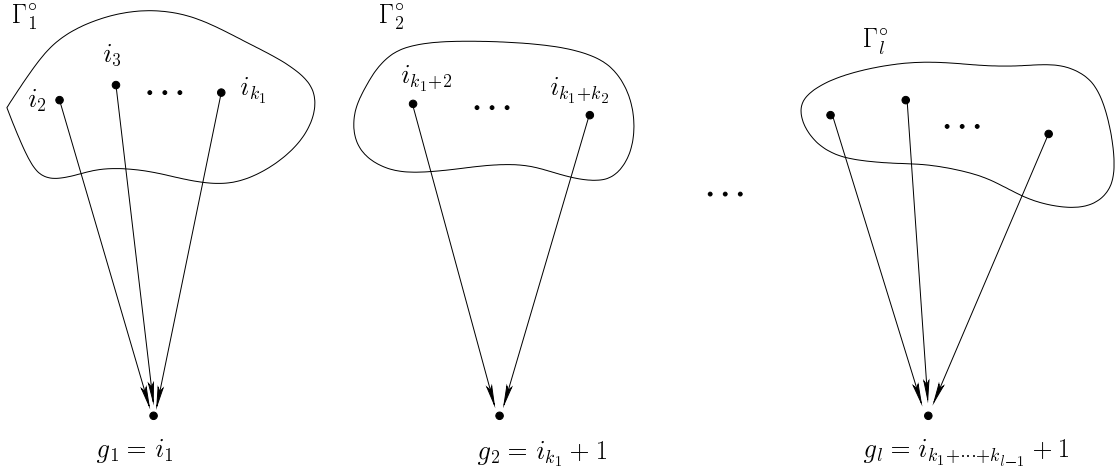}

\caption{Correspondence\label{cap:or-graph} between maps $\pi _{k_{1},\ldots ,k_{l}}$
and oriented graphs}
\end{figure}

Let us enumerate the particles and let $x_{1},\ldots ,x_{N}$ be their
coordinates. Let us fix a natural number~$k\geq 2$ and integers
$k_{1}\geq 2$, $\ldots $, $k_{l}\geq 2$ such that $k_{1}+\cdots +k_{l}=k$.
We call the sequenced collection $(k_{1},\ldots ,k_{l})$ a \emph{signature}
and will keep it fixed throughout the present paper. Let $\mathcal{I}$
be a set of all sequenced collections $(i_{1},\ldots ,i_{k})$, consisting
of $k$ different elements of the index set $\nN :=\left\{ 1,\ldots ,N\right\} $.
On the set $\mathcal{I}$ we define a map $\pi _{k_{1},\ldots ,k_{l}}:\, \, (i_{1},\ldots ,i_{k})\mapsto \left(\Gamma _{1},\ldots ,\Gamma _{l}\right),$
where $\Gamma _{j}=(g_{j},\Gamma _{j}^{\circ })$, \begin{eqnarray*}
g_{1}=i_{1},\quad \quad \Gamma _{1}^{\circ } & = & (i_{2},\ldots ,i_{k_{1}}),\\
g_{2}=i_{k_{1}+1},\quad \quad \Gamma _{2}^{\circ } & = & (i_{k_{1}+2},\ldots ,i_{k_{1}+k_{2}}),\\
\cdots  &  & \\
g_{l}=i_{k_{1}+\cdots +k_{l-1}+1},\quad \quad \Gamma _{l}^{\circ } & = & \left(i_{k_{1}+\cdots +k_{l-1}+1},\ldots ,i_{k_{1}+\cdots +k_{l}}\right).
\end{eqnarray*}
It is useful to associate with the map $\pi _{k_{1},\ldots ,k_{l}}$
an oriented graph as shown on the Figure~\ref{cap:or-graph}. %

Consider also a map $J_{k_{1},\ldots ,k_{l}}^{\ioik }:\, \, \RR ^{N}\rightarrow \RR ^{N}$,\[
J_{k_{1},\ldots ,k_{l}}^{\ioik }:\quad x=\left(x_{1},\ldots ,x_{N}\right)\mapsto y=\left(y_{1},\ldots ,y_{N}\right),\]
which is defined as follows:\begin{equation}
y_{m}=\left\{ \begin{array}{rl}
 x_{m}, & \esli m\notin \ioik ,\\
 x_{g_{j}}, & \esli m\in \ioik ,\, m\in \Gamma _{j}\, .\end{array}\right.\label{eq:yJx}\end{equation}
It is natural to call this map a \emph{synchronization} of the collection
of particles $x_{i_{1}},\ldots ,x_{i_{k}}$, \emph{corresponding}
to the signature $(k_{1},\ldots ,k_{l})$. It is evident that the
configuration $J_{k_{1},\ldots ,k_{l}}^{\ioik }x$ has at least $k_{1}$
particles with coordinates that are equal to $x_{g_{1}}$, \ldots,
at least $k_{l}$ particles with coordinates that are equal to $x_{g_{l}}$.

\paragraph*{\punkt Multidimensional Markov process with synchronization. }

We define a continuous time homogeneous Markov process $x(t)=\left(x_{1}(t),\ldots ,x_{N}(t)\right)$,
$t\geq 0$, on the state space $\RR ^{N}$ by means of the following
formal generator $L=\Lf +\Ls $. \emph{Free dynamics} $\Lf $ corresponds
to independent random walks and is chosen as \[
\left(\Lf f\right)(x)=\al \sum _{i=1}^{N}\int \left(f\left(x+ze_{i}\right)-f\left(x\right)\right)\rho (dz),\qquad z\in \RR ,\quad e_{i}=\left(0,\ldots ,1_{i},\ldots ,0\right),\]
i.e., independently of the other particles each particle~$i$ waits
an exponential time with parameter $\al >0$ and performs a jump of
the form $x_{i}\rightarrow x_{i}+z$, where $z$ is distributed according
to the law $\rho (dz)$, which is common for all particles. We assume
that the distribution $\rho $ has a compact support and is nontrivial
in the following sense: $\sdv :=\int x^{2}\rho (dx)>0$. Denote also
$a:=\int x\rho (dx).$

\emph{Synchronizing interaction $\Ls $} has the form

\begin{equation}
\Ls f=\frac{\de }{N^{[k]}}\sum _{(i_{1},\ldots ,i_{k})\in \mathcal{I}}\Lsii f,\label{eq:Ls-ioik}\end{equation}
where $N^{[k]}:=N(N-1)\cdots (N-k+1)$ and\begin{equation}
\left(\Lsii f\right)(x)=f\left(J_{k_{1},\ldots ,k_{l}}^{\ioik }x\right)-f(x),\qquad x=(x_{1},\ldots ,x_{N}),\label{eq:Lsf}\end{equation}
for bounded continuous functions $f$. In other words, independently
from the free dynamics at epochs of a Poissonian flow with parameter~~$\de >0$
we choose with an equal probability one element $\ioik $ from the
set $\mathcal{I}$, and synchronize the particle configuration $x_{1},\ldots ,x_{N}$
according to the map $J_{k_{1},\ldots ,k_{l}}^{\ioik }$. Hence, the
interaction, corresponding to the given choice of $\ioik $, consists
in instantaneous and simultaneous transfer of particles of each group
$x_{k_{1}+\cdots +k_{j-1}+1},\ldots ,x_{k_{1}+\cdots +k_{j}}$ to
the point with the coordinate $x_{g_{j}}=x_{k_{1}+\cdots +k_{j-1}+1}$,
$j=\overline{1,l}$. 

Since in the sum~(\ref{eq:Ls-ioik}) all summands have the same weight,
the system under consideration belongs to the class of models with
mean-field interaction.

\paragraph*{\punkt Main results.}

Consider a function $V:\, \RR ^{N}\rightarrow \RR _{+}$\[
V(x):=\frac{1}{N(N-1)}\sum _{m<n}\left(x_{m}-x_{n}\right)^{2}.\]
It is easy to see that this function coincides with the \emph{empiric
(sample) variance}~$S^{2}$, where\[
S^{2}:=\frac{1}{N-1}\sum _{m=1}^{N}\left(x_{m}-M(x)\right)^{2},\qquad M(x):=\frac{1}{N}\sum _{m=1}^{N}x_{m}\, .\]

\begin{lemma}\label{l-LM}

$\left(LM\right)(x)=\al a\, $ for all $x=(x_{1},\ldots ,x_{N})$.

\end{lemma}

\begin{lemma}\label{l-LV}

If $N\rightarrow \infty $ then we have $\ds LV=\al \sdv -\frac{\de \cdv }{N(N-1)}\, V$
for some $\cdv >0$.

\end{lemma}

It follows from the proof that $\cdv =\ds \sum _{j=1}^{l}k_{j}^{2}-k$.

We are interested in the following means: $\mux _{N}(t):=\sE M(x(t))$
and $\erx _{N}(t):=\sE V(x(t))$. In the below statements we always
assume that $N\rightarrow \infty $.

\begin{theorem}\label{t-mu}

For any $t>0\, \, $ $\ds \lim _{N}\frac{\mu _{N}(t)-\mu _{N}(0)}{t}=\al a$.
Moreover, for any function $t(N)\rightarrow \infty $ the following
convergence holds $\ds \frac{\mu _{N}(t(N))-\mu _{N}(0)}{t(N)}\rightarrow \al a$
.

\end{theorem}

\begin{theorem}\label{t-R}

Assume that $\sup _{N}\left|\erx _{N}(0)\right|<\infty $. There exist
three time scales $t=t(N)$, where $R_{N}(t(N))$ has completely different
asymptotic behavior.

\begin{itemize}
\item If $\ds \frac{t(N)}{N^{2}}\rightarrow 0$, then $R_{N}(t(N))\sim \, \al \sdv \, t(N)$.
\item If $t(N)=cN^{2}$, then $\ds R_{N}(t(N))\sim \al \sdv (\de \cdv )^{-1}(1-\exp \left(-\de \cdv c\right))N^{2}$.
\item If $\ds \frac{t(N)}{N^{2}}\rightarrow \infty $, then $\ds R_{N}(t(N))\sim \al \sdv (\de \cdv )^{-1}N^{2}$.
\end{itemize}
\end{theorem} In~\cite{malysh-Manita-tvp} similar consecutive stages
were called correspondingly a phase of initial desynchronization,
a phase of critical slowdown of desynchronization and a phase of final
stabilization. We see that the first phase does not contribute to
the asymptotic behavior.

\begin{theorem}\label{t-sled}Assume that $\erx _{N}(0)\rightarrow \infty $
as $N\rightarrow \infty $. Then the additional condition $\ds \frac{\erx _{N}(0)}{t(N)}\rightarrow 0$
is sufficient to ensure the validity of the corresponding statements
of Theorem~\ref{t-R}.

\end{theorem}

\paragraph*{\punkt Proof of Lemma~\ref{l-LM}. }

A straightforward calculation shows that $\left(\Lf M\right)(x)\equiv \al a$,
hence we need only to prove that $\left(\Ls M\right)(x)\equiv 0$.
Let us show that this fact follows from the symmetry of the synchronizing
interaction. Indeed, if the signature of interaction $(k_{1},\ldots ,k_{l})$
is given we can fix some sets\[
B_{j}\in \nN ,\quad \left|B_{j}\right|=k_{j},\quad (j=\overline{1,l}),\quad \quad B_{j_{1}}\cap B_{j_{2}}=\varnothing \quad (j_{1}\not =j_{2})\]
and define $\mathcal{I}_{k_{1},\ldots ,k_{l}}^{B_{1},\ldots ,B_{l}}:=\left\{ \ioik :\, \, \left\{ g_{j}\right\} \cup \Gamma _{j}^{\circ }=B_{j}\, \, \, \forall j=\overline{1,l}\right\} ,$
where $g_{j}$ and $\Gamma _{j}^{\circ }$ are determined by the map
$\pi _{k_{1},\ldots ,k_{l}}$ applied to $\ioik $. To finish the
proof it is sufficient to show that\begin{equation}
\sum _{\ioik \in \mathcal{I}_{k_{1},\ldots ,k_{l}}^{B_{1},\ldots ,B_{l}}}\left(M\circ J_{k_{1},\ldots ,k_{l}}^{\ioik }-M\right)=0,\label{eq:B-B-0}\end{equation}
since $\sum \limits _{(i_{1},\ldots ,i_{k})\in \mathcal{I}}=\sum \limits _{B_{1},\dots ,B_{l}}\, \sum \limits _{\ioik \in \mathcal{I}_{k_{1},\ldots ,k_{l}}^{B_{1},\ldots ,B_{l}}}$.
It is easy to see that \begin{equation}
M\left(J_{k_{1},\ldots ,k_{l}}^{\ioik }x\right)-M(x)=\frac{1}{N}\sum _{j=1}^{l}\sum _{h\in \Gamma _{j}^{\circ }}\left(x_{g_{j}}-x_{h}\right).\label{eq:M-J-M}\end{equation}
Substituting the r.h.s.~of~(\ref{eq:M-J-M}) into the sum~(\ref{eq:B-B-0}),
we see that, since in~(\ref{eq:B-B-0}) $g_{j}$ covers the whole
set $B_{j}$, for any pair of indices $u,v\in B_{j}$ we have exactly
one difference $(x_{u}-x_{v})$, when $g_{j}=u$, and exactly one
difference $(x_{v}-x_{u})$, when $g_{j}=v$. Hence the sum~(\ref{eq:B-B-0})
is equal to zero.

\paragraph*{\punkt Proof of Lemma~\ref{l-LV}. }

Define a function $f_{m,n}:\, \RR ^{N}\rightarrow \RR $ by the formula
\[
f_{m,n}(x)=\ds \frac{1}{N(N-1)}(x_{n}-x_{m})^{2}.\]
So $V(x)=\ds \sum _{m<n}f_{m,n}(x)$. It is straightforward to check
that $\Lf f_{m,n}(x)=\Bigl (N(N-1)\Bigr )^{-1}2\al \sdv $ and hence
$\left(\Lf V\right)(x)\equiv \al \sdv $.

We have \begin{equation}
\Ls V=\de \sum _{m'<n'}\frac{1}{N^{[k]}}\sum _{\ioik \in \mathcal{I}}\Lsii f_{m',n'}\, .\label{eq:LVLf}\end{equation}
Consider a summand $\Lsii f_{m',n'}$. By~(\ref{eq:yJx}) the map
$J_{k_{1},\ldots ,k_{l}}^{\ioik }$ transfers a particle with index~$m'$
to the point having the coordinate $x_{m}$, where\begin{equation}
m=\left\{ \begin{array}{rl}
 m', & \esli m'\notin \ioik ,\\
 g_{j}, & \esli m'\in \ioik ,\, m'\in \Gamma _{j}\end{array}\right.\label{eq:m-sh-m}\end{equation}

Consequently, $\Lsii f_{m',n'}=f_{m,n}-f_{m',n'}$ with some $m$
and $n$, which are not necessarily different from $m'$ and $n'$.
Hence, there exist such $a_{nm}(N)\in \RR $ that\begin{equation}
\left(\Ls V\right)(x)=\de \sum _{m<n}a_{mn}(N)f_{m,n}(x).\label{eq:Laf}\end{equation}
Our goal is to show that the coefficients $a_{mn}(N)$ do not depend
on $m$ and $n$ and that\begin{equation}
a_{mn}(N)=\, -\frac{\cdv }{N(N-1)}\, \label{eq:anmC2}\end{equation}
for some constant $\cdv >0$. Let us fix some pair $\left\{ m,n\right\} $
and calculate~$a_{mn}(N)$. When we choose a collection $\ioik $
for the synchronization, in that way we choose $g_{1},\ldots ,g_{l}$
and $\Gamma _{1}^{\circ },\ldots ,\Gamma _{l}^{\circ }$. Denote $G=\left\{ g_{1},\ldots ,g_{l}\right\} $.
If the pair $\{m,n\}$ is given then the set~$\mathcal{I}$ can be
divided into four non-intersecting parts as follows: $\mathcal{I}=\ds \bigcup _{w=0}^{3}\mathcal{I}_{w}^{m,n}$,
where\begin{eqnarray*}
\mathcal{I}_{0}^{m,n}: & = & \left\{ \ioik :\, \left\{ m,n\right\} \cap \left(\cup _{j}\Gamma _{j}\right)=\varnothing \right\} \, ,\\
\mathcal{I}_{1}^{m,n}: & = & \left\{ \ioik :\, \left\{ m,n\right\} \cap \left(\cup _{j}\Gamma _{j}^{\circ }\right)\not =\varnothing \right\} \, ,\\
\mathcal{I}_{2}^{m,n}: & = & \left\{ \ioik :\, \left\{ m,n\right\} \subset G\right\} \, ,\\
\mathcal{I}_{3}^{m,n}: & = & \left\{ \ioik :\, \left|\left\{ m,n\right\} \cap G\right|=1\right\} \backslash \mathcal{I}_{1}^{m,n}\, .
\end{eqnarray*}
In each sum $\sum ^{(w)}:=\ds \, \frac{1}{N^{[k]}}\, \sum _{\ioik \in \mathcal{I}_{w}^{m,n}}\, \sum _{m'<n'}\left(f_{m',n'}\circ J_{k_{1},\ldots ,k_{l}}^{\ioik }-f_{m',n'}\right)$
we pick out only summands that contain the function $f_{m,n}$ and
a coefficient in front of it will be denoted by~$a_{mn}^{(w)}(N)$.
By representations~(\ref{eq:LVLf}) and~(\ref{eq:Laf}) we can write
$a_{nm}(N)=\sum \limits _{w=0}^{3}a_{mn}^{(w)}(N).$

0) If $\ioik \in \mathcal{I}_{0}^{m,n}$ then the particles with indices~$m$
and $n$ are fixed points of the map $J_{k_{1},\ldots ,k_{l}}^{\ioik }$,
hence $a_{mn}^{(0)}(N)=0$.

I) If $\ioik \in \mathcal{I}_{1}^{m,n}$, then the summand $f_{m,n}$
can be presented in the sum~$\sum ^{(1)}$ only with the sign {}``$\, -\, $''
and only in the case when $(m',n')=(m,n)$. A total number of collections
$\ioik $, that belong to the set $\mathcal{I}_{1}^{m,n}$, is equal
to $\left|\mathcal{I}_{1}^{m,n}\right|=N^{[k]}-\left(N-2\right)^{[k-l]}\left(N-(k-l)\right)^{[l]}$.
Therefore, \[
a_{mn}^{(1)}(N)=\left(-1\right)\frac{N^{[k]}-\left(N-2\right)^{[k-l]}\left(N-(k-l)\right)^{[l]}}{N^{[k]}}\, .\]

II) For definiteness let us fix some $i,j\in \left\{ 1,\ldots ,l\right\} $
and consider a subsum of $\sum ^{(2)}$ taken over such subset of
$\mathcal{I}_{2}^{m,n}$ that $m=g_{i}$, $n=g_{j}$. Under the action
of the map $J_{k_{1},\ldots ,k_{l}}^{\ioik }$ each of function $f_{uv}$,
where $u\in \Gamma _{i}$, $v\in \Gamma _{j}$, turns into the function
$f_{mn}$. Thus in the above subsum we find $k_{i}k_{j}$ summands
$f_{mn}$ with the sign {}``$\, +\, $'' and only one summand (corresponding
to the case $(m',n')=(m,n)\, $) with the sign {}``$\, -\, $''.
Hence, \[
a_{mn}^{(2)}(N)=\sum _{i,j=1,i\not =j}^{l}\left(k_{i}k_{j}-1\right)\frac{\left(N-2\right)^{[k-2]}}{N^{[k]}}\, .\]

III) Fix $i\in \left\{ 1,\ldots ,l\right\} $ and consider in $\sum ^{(3)}$
a subsum taken over such~$\ioik $ that $m=g_{i}\in G$, $n\notin G$.
A number of collections $\ioik $, satisfying to this assumption,
is equal to $(N-2)^{[k-1]}$. Under this assumption each of $k_{i}$
functions $f_{un}$, where $u\in \Gamma _{i}$, turns into $f_{mn}$.
Changing the roles of $m$ and $n$, we conclude that\[
a_{mn}^{(3)}(N)=2\sum _{i=1}^{l}\left(k_{i}-1\right)\frac{(N-2)^{[k-1]}}{N^{[k]}}\, .\]

As is easy to see the values of $a_{mn}^{(1)}(N)$, $a_{mn}^{(2)}(N)$
and $a_{mn}^{(3)}(N)$ do not depend on $m,n$. After some calculations
we obtain the following formulae \begin{eqnarray*}
a_{mn}^{(1)}(N) & = & -\frac{2(k-l)(N-k)+(k-l)(k+l-1)}{N(N-1)}\, ,\\
a_{mn}^{(2)}(N) & = & \frac{k^{2}-l^{2}+l-\sum _{j}k_{j}^{2}}{N(N-1)}\, ,\\
a_{mn}^{(3)}(N) & = & \frac{2(k-l)(N-k)}{N(N-1)}\, .
\end{eqnarray*}
Summing these values together, we get\[
a_{nm}^{(1)}(N)+a_{nm}^{(2)}(N)+a_{nm}^{(3)}(N)=-\frac{\sum _{j}k_{j}^{2}-k}{N(N-1)}\, ,\]
and the statement~(\ref{eq:anmC2}) is proved with $\cdv =\sum _{j}k_{j}^{2}-k>0$.

\paragraph*{\punkt Proofs of theorems.}

Our method is similar to the approach which was proposed in~\cite{malysh-Manita-tvp}.
It is based on an embedded Markov chain $\zeta _{N}(n,\om ):=x(\tau _{n},\om ),$
$n=0,1,\ldots \, $. The chain $\zeta _{N}(n)$ is defined on the
same probability space $\left(\Omega ,\mathcal{F},\mathsf{P}\right)$,
evolves in discrete time and take its values in the set~$\RR ^{N}$.
A sequence $\tau _{1}(\om )\leq \tau _{2}(\om )\leq \cdots $ consists
of time moments at which particles make jumps, $\tau _{0}(\om )\equiv 0$.
Evidently, $\left\{ \tau _{n+1}-\tau _{n}\right\} _{n=0}^{\infty }$
is a sequence of i.i.d~random variables having exponential distribution
with the mean value~$\gn =\left(\al N+\de \right)^{-1}$. In the
discrete Markov chain~$\zeta _{N}(n)$ one-step transitions have
the following form: with probability $\al \gn $ particle~$i$ makes
a jump $x_{i}\rightarrow x_{i}+z$, where $z$ is distributed according
to the law~$\rho (dz)$, $i=\overline{1,N}$, or with probability
$\delta \gn $ a synchronization take place, namely, a collection
of indices~$\ioik $ is chosen randomly and then the particle configuration
$x$ is transformed into~$J_{k_{1},\ldots ,k_{l}}^{\ioik }x$. By
the law of large numbers $\tau _{n}\sim \gn n$ as $n\rightarrow \infty $,
therefore, an asymptotic behavior of the particle system~$x(t)$,
$t\geq 0$, can be reduced to asymptotic properties of the embedded
Markov chain~$\zeta _{N}(n)$. We are interested in the following
sequences: $\srM (n):=\sE M\left(x(\tau _{n})\right)$ and $\dsp (n):=\sE V(x(\tau _{n}))$.
Straightforward calculations together with Lemmas~\ref{l-LM} and~\ref{l-LV}
show that \begin{eqnarray*}
\sE \left(M(\zeta _{N}(n+1))\, \, |\, \, x(t),\, t\leq \tau _{n}\right) & = & M(\zeta _{N}(n))+\gn \al a,\\
\sE \left(V(\zeta _{N}(n+1))\, \, |\, \, x(t),\, t\leq \tau _{n}\right) & = & V(\zeta _{N}(n))+\gn \left(\al \sdv -\de \frac{\cdv }{N(N-1)}V(\zeta _{N}(n))\right).
\end{eqnarray*}
Taking expected values of the both parts in the above equations we
come to recurrent relations\begin{eqnarray}
\srM (n+1) & = & \srM (n)+\gn \al a,\label{eq:sr-lin}\\
\dsp (n+1) & = & \dsp (n)\left(1-\gn \de \frac{\cdv }{N(N-1)}\right)+\gn \al \sdv \label{eq:disp-rec}
\end{eqnarray}
From~(\ref{eq:sr-lin}) the statement of Theorem~\ref{t-mu} easily
follows. By iterating the equation~(\ref{eq:disp-rec}), we get\begin{equation}
\dsp (n)=d(0)\left(1-\gn \de \frac{\cdv }{N(N-1)}\right)^{n}+\left(\gn \al \sdv \right)\sum _{j=1}^{n-1}\left(1-\gn \de \frac{\cdv }{N(N-1)}\right)^{j}\, .\label{eq:dn-d0}\end{equation}
Substituting $n=\gn ^{-1}t(N)$ in~(\ref{eq:dn-d0}) and letting
$N$ go to the infinity, we come to the conclusion of Theorem~\ref{t-R}.

Theorem~\ref{t-sled} can be derived from~(\ref{eq:dn-d0}) by a
straightforward calculation performed separately for each time stage.

\paragraph*{\punkt Possible generalizations and perspectives. }

Undoubtedly the results obtained in this paper remain also true, if
the synchronization $\Ls $ is taken in the form of a symmetrized
polynomial interaction of any order, i.e., in the case when $\Ls $
is a finite sum (taken over index $k\geq 2$) of $k$-particle symmetrized
interactions. Apparently, the methods of our paper can be adapted
also for other classes of the free evolutions~$\Lf $, such as independent
Brownian motions.

The stochastic particle system considered here plays a basic role
for future study of asymptotic behavior of general many-component
stochastic systems with synchronization. The statement of Theorem~\ref{t-R}
put forward a hypothesis that, seemingly, the result of the paper~\cite{malysh-Manita-tvp}
on the existence of three phases of collective behavior remains true
for a wide class of large particle systems with synchronization and,
hence, it is not really caused by the specific nature of the model~\cite{malysh-Manita-tvp}.
As it is seen from the present paper, key elements in the future proofs
of such results should be some analogues of Lemma~\ref{l-LV}. On
this way we expect to have difficulties with anisotropic synchronizations
which are interesting in a number of important applications. Let us
remark that in papers~\cite{Man-Sche,Malyshkin}, where some examples
of anisotropic interactions were considered, the behavior of particle
systems was studied only on a so called hydrodynamic scale while in
the present paper we consider all possible time scales $t(N)$.

Unfortunately, an important class of cascade synchronizations~\cite{Man-Sim1,Man-Sim}
can not be considered in the framework of the present approach.

~
\end{document}